\documentclass[12pt]{amsart}

\newtheorem{theorem}{Theorem}[section]
\newtheorem{lemma}[theorem]{Lemma}

\begin{document}
\baselineskip=15.5pt

\newcommand{\wh}[1]{\widehat{#1}}

\newcommand{\pardeg}{{\rm deg_{par}}}

\newcommand{\Ocal}{{\cal O}}

\newcommand{\cdop}{{\mathbb C}}
\newcommand{\ndop}{{\mathbb N}}
\newcommand{\pdop}{{\mathbb P}}
\newcommand{\zdop}{{\mathbb Z}}

\title{Parabolic Raynaud bundles}

\author[I. Biswas]{Indranil Biswas}

\address{School of Mathematics, Tata Institute of Fundamental
Research, Homi Bhabha Road, Bombay 400005, India}

\email{indranil@math.tifr.res.in}

\author[G. Hein]{Georg Hein}

\address{Universit\"at Duisburg-Essen, Fachbereich
Mathematik, 45117 Essen, Germany}

\email{georg.hein@uni-due.de}

\subjclass[2000]{14F05, 14H60}

\keywords{Parabolic bundle, semistability, Raynaud bundle}

\date{September 13, 2007}

\begin{abstract}

Let $X$ be an irreducible smooth projective curve defined over
complex numbers, $S=\{p_1, p_2,\ldots,p_n\} \subset X$ a finite
set of  closed points and $N\,  \geq\, 2$ a fixed integer. For any
pair $(r,d) \in \ndop \times \frac{1}{N} \zdop$, there exists a
parabolic vector bundle $R_{r,d,*}$ on $X$, with parabolic structure
over $S$ and all parabolic weights in $\frac{1}{N} \zdop$, that
has the following property:
Take any parabolic vector bundle $E_*$ of rank $r$ on $X$ whose
parabolic points are contained in $S$,
all the parabolic weights are in $\frac{1}{N}\zdop$
and the parabolic degree is $d$.
Then $E_*$ is parabolic semistable if and only if there is
no nonzero parabolic homomorphism from $R_{r,d,*}$ to $E_*$.

\end{abstract}

\maketitle

\section{Introduction}

Let $X$ be an irreducible smooth projective curve defined over
$\mathbb C$. Faltings showed that a vector bundle $E$ over $X$
is semistable if and only if there exists a vector bundle $F$
on $X$
such that $H^i(X, {\mathcal H}om(F,E)) =0$ for $i=0,1$ \cite{Fal}.
This result was developed further by Popa (see \cite{Pop})
to obtain estimates for the rank of $F$ which depend
only on the rank of $E$. 
Given $r$ and $d$, a vector bundle $R_{r,d}$ on $X$ is called
a \textit{Raynaud bundle} if the following holds: A vector bundle
$E$ on $X$ of rank $r$ and degree $d$ is semistable if and only
if $\text{Hom}(R_{r,d}, E) =0$ \cite{Hei}.
We note that 
these are called Raynaud bundles because the vector bundles
considered by Raynaud in \cite{Ray} are Raynaud bundles for
$r\, =\,2$ and $d\, =\, 2(\text{genus}(X)-1)$
(see Theorem 1.1 in \cite{Hei}). 

In \cite{Hei},
building on Popa's result it was shown that Raynaud bundles exist.
In \cite{Bi3}, an analog of Faltings' semistability criterion 
is given for parabolic vector bundles. Our aim here is to prove
the existence of the Raynaud bundles in the parabolic context.

Given two parabolic vector bundles $E_*$ and $F_*$ on $X$, the
global parabolic homomorphisms from $F_*$ to $E_*$ will be
denoted by $\text{Hom}_{\rm par}(F_*, E_*)$ (see
\cite[p. 212, Definition 1.5(II)]{MS} for the definition of
a parabolic homomorphism).

We prove the following theorem:

\begin{theorem}\label{THE-1}
Let $S=\{p_1, p_2,\ldots,p_n\} \subset X$ a finite set of closed
points  of $X$, and $N\, \geq\, 2$ a fixed integer. For any pair
$(r,d) \in \ndop \times \frac{1}{N} \zdop$, there exists a
parabolic vector bundle $R_{r,d,*}$ on $X$, with parabolic
structure over $S$ and all parabolic weights in
$\frac{1}{N} \zdop$, that has the following property:
Take any parabolic vector bundle $E_*$ on $X$ such that
\\
\begin{tabular}{lp{13cm}}
(1) & the parabolic points of $E_*$ are contained in $S$,\\
(2) & all the parabolic weights of $E_*$ are in $\frac{1}{N} \zdop$, 
and\\
(3) & the rank of $E_*$ is $r$, and the parabolic degree of $E_*$ is
$d$.\\
\end{tabular}\\
Then $E_*$ is parabolic semistable $\iff \,
{\rm Hom}_{\rm par}(R_{r,d,*}, E_*)\,=\, 0$.
\end{theorem}

\section{Preliminaries}\label{sec2}

\subsection{The equivalence of equivariant bundles and parabolic
bundles}\label{suse-1}

We will recall a correspondence between parabolic vector
bundles and equivariant vector bundles which will be used
in the proof of Theorem \ref{THE-1}. 

We assume that at least one of the following two conditions hold:
\begin{itemize}
\item $\text{genus}(X)\, \geq\, 1$
\item $|S|\, \not=\, 1$.
\end{itemize}

Fix $S$ and $d$ as in Theorem \ref{THE-1}. Fix 
a Galois algebraic covering
\begin{equation}\label{f}
f\,:\, Y\, \longrightarrow\, X
\end{equation}
such that
\begin{itemize}
\item $f$ is ramified exactly over $S$, and

\item the ramification index of each point in $f^{-1}(S)$ is
$N-1$.
\end{itemize}
See \cite[p. 26, Proposition 1.2.12]{Na} for the
existence of $f$ satisfying these conditions. We note that
the assumption $(\text{genus}(X)\, , |S|)\, \not=\,
(0\, ,1)$ is needed for the existence of $f$.

Let
$$
\Gamma\, :=\, \text{Gal}(f)
$$
be the Galois group for the covering $f$. A
$\Gamma$--\textit{linearized vector bundle} on $Y$ is an algebraic
vector bundle $E$ equipped with a lift of the action of
$\Gamma$ as vector bundle automorphisms. This means that $\Gamma$
acts on the total space of $E$ as algebraic automorphisms, and
the action of each $\gamma\, \in\, \Gamma$ on $E$ is an
isomorphism of the vector bundle $E$ with $(\gamma^{-1})^*E$.

In \cite{Bis1}, a natural bijective correspondence between the
following two classes was established:
\begin{enumerate}
\item the $\Gamma$--linearized vector bundles $W$ on $Y$, and 

\item the parabolic vector bundles $E_*$ over $X$ for which the
parabolic divisor in contained in $S$ and all the
parabolic weights are in $\frac{1}{N} \zdop$.
\end{enumerate}

This bijective correspondence takes the usual tensor
product (respectively, direct sum) of 
$\Gamma$--linearized vector bundles to the
tensor product (respectively, direct sum) of the
corresponding parabolic vector bundles. Similarly,
for any $\Gamma$--linearized vector bundle $W$,
this bijective correspondence takes $W^*$ to the parabolic
dual of the parabolic vector bundle corresponding to $W$.
(See \cite{Bi2}, \cite{Yo} for the above mentioned operations
on parabolic vector bundles.)

For any parabolic vector bundle $E_*$ of the above type, if
$\wh E$ is the corresponding $\Gamma$--linearized vector
bundle, then
\begin{equation}\label{eq1}
\deg(\wh E) = | \Gamma| \cdot \pardeg(E_*)\, ,
\end{equation}
where $\pardeg(E_*)$ is the parabolic degree of $E_*$.
\cite[p. 318, (3.12)]{Bis1}. Furthermore,
\begin{equation}\label{eq2}
E_* \mbox{ is parabolic semistable} \iff \wh E \mbox{ is semistable}
\end{equation}
(see \cite[p. 318, Lemma 3.13]{Bis1}).

\subsection{Raynaud bundles on smooth algebraic curves}\label{suse-2}
In \cite{Hei} an irreducible smooth projective curve $Y$
is considered. It is shown that for any pair
of integers $(r,d)$, there
exists a vector bundle $R_{r,d}$ on $Y$ (which is called a
Raynaud bundle) with the following property:

For a vector bundle $E$ of rank $r$ and degree $d$ on $Y$,
\begin{equation}\label{eq3}
E \mbox{ is semistable} \iff {\rm Hom}(R_{r,d},E)
\,=\,0
\end{equation}
(see (i) $\iff$ (v) in Theorem 2.12 of \cite{Hei}).

This way we obtain a short cohomological criterion which
enables us to check semistability.
We can compute the rank and degree of $R_{r,d}$
in terms of the genus $g_Y$ of $Y$ and the integers $r$ and
$d$ (see Proposition 2.2 and Corollary 3.4 in \cite{Hei}).

\section{Proof of Theorem \ref{THE-1}}

As in Section \ref{suse-1}, we will assume that at
least one of the following two conditions hold:
\begin{itemize}
\item $\text{genus}(X)\, \geq\, 1$
\item $|S|\, \not=\, 1$.
\end{itemize}
The remaining case where $(\text{genus}(X)\, , |S|)\, =\,
(0\, ,1)$ will be treated separately.

Let $E_*$ be a parabolic vector
bundle on $X$ satisfying conditions (1)--(3) of
Theorem \ref{THE-1}. Denoting by $\wh E$ the associated
$\Gamma$--linearized vector bundle on
$Y$ in (\ref{f}), we have that $\wh E$ is of rank $r$ and degree
$d| \Gamma|$, and the parabolic semistability of $E_*$ is
equivalent to the semistability of $\wh E$ (see (\ref{eq2})).

Using the result (\ref{eq3}) together with (\ref{eq1})
and (\ref{eq2})
we obtain that
\begin{equation}\label{e0}
E_* \mbox{ is parabolic semistable }\, \iff \,
{\rm Hom}(R_{r,d| \Gamma|},\wh E)\,=\,0\, .
\end{equation}

{\bf Step 1:}
Set now
\begin{equation}\label{e-1}
{\widetilde R}_{r,d} = \bigoplus_{\gamma \in \Gamma}
\gamma^*R_{r,d}
\end{equation}
to be the direct sum. Since $\wh E$ is $\Gamma$--linearized,
$$
{\rm Hom}(\widetilde R_{r,d| \Gamma|},\wh E)\, \cong\,
{\rm Hom}(R_{r,d| \Gamma|},\wh E)^{\oplus |\Gamma|}\, .
$$
Therefore, 
$$
{\rm Hom}(\widetilde R_{r,d| \Gamma|},\wh E)\,=\,0
\iff
{\rm Hom}(R_{r,d| \Gamma|},\wh E)\,=\,0\, .
$$
Combining this with (\ref{e0}) we conclude that
\begin{equation}\label{eq4}
E_* \mbox{ is parabolic semistable }\,
\iff \,
{\rm Hom}(\widetilde R_{r,d| \Gamma|},\wh E)\,=\,0\, .
\end{equation}

{\bf Step 2:}
The vector bundle ${\widetilde R}_{r,d}$ in (\ref{e-1})
admits a canonical $\Gamma$--linearization.
Let $R'_{r,d,*}$ be the parabolic vector bundle over
$X$ corresponding to this $\Gamma$--linearized vector bundle
${\widetilde R}_{r,d}$.

Consider the trivial vector bundle over $Y$
\begin{equation}\label{hW}
\widehat{W}\, :=\, {\mathcal O}_Y\otimes_{\mathbb C}
{\mathbb C}(\Gamma)\, ,
\end{equation}
where ${\mathbb C}(\Gamma)$ is the group algebra of $\Gamma$.
The action of $\Gamma$ on $Y$ lifts to an action of $\Gamma$
on ${\mathcal O}_Y$. The natural action of $\Gamma$
on ${\mathbb C}(\Gamma)$ and the action of $\Gamma$
on ${\mathcal O}_Y$ together define a 
$\Gamma$--linearization on the vector bundle $\widehat{W}$
in (\ref{hW}).

Let $W_*$ denote the parabolic vector bundle over $X$ corresponding
to the $\Gamma$--linearized vector bundle $\widehat{W}$.

\begin{lemma}\label{lem1}
Let $F_*$ be a parabolic vector bundle over $X$ satisfying
condition (1) and condition (2) in Theorem \ref{THE-1}. Let
$\widehat{F}$ be the $\Gamma$--linearized vector bundle on $Y$
corresponding to $F_*$. Then
$$
{\rm Hom}({\widetilde R}_{r,d| \Gamma|},\wh F)\, =\,
{\rm Hom}_{\rm par}(R'_{r,d,*}\otimes W_*, F_*)\, ,
$$
where $R'_{r,d,*}\otimes W_*$ is the parabolic tensor product of
the parabolic vector bundles $R'_{r,d,*}$ and $W_*$ constructed
above.
\end{lemma}

\noindent
\textbf{Proof.} 
The parabolic vector bundle $W_*$ is constructed as follows.
Consider the direct image
\begin{equation}\label{W}
W\, =\, f_* {\mathcal O}_Y\, ,
\end{equation}
where $f$ is the covering map in (\ref{f}).
We have a filtration of subsheaves
\[
W_1\,\subset\,\cdots\, \subset\, W_i\,
\subset \,\cdots\, \subset\, W_{N-1}\,\subset\, W_N \, =\, W\, ,
\]
where $W_j\, :=\, f_* {\mathcal O}_Y(-(N-j)f^{-1}(S)_{\text{red}})$.
For any point $x\, \in\, S$, let
\begin{equation}\label{eq4-1}
0\, \subset\,
W^1_{x} \,\subset\,\cdots\, \subset\, W^j_{x}\,\subset\,
W^{j+1}_{x} \,\subset\,\cdots\, \subset\, W^{N-1}_{x}\,
\subset\,W^N_{x} \, =\, W_{x}
\end{equation}
be the filtration of subspaces given by the above filtration
of subsheaves. The dimension of each successive
quotient in (\ref{eq4-1}) is $|\Gamma|/N$. 

The vector bundle underlying the parabolic vector bundle
$W_*$ is $W$ (defined in (\ref{W})), its parabolic
divisor is $S$, its quasiparabolic filtration on each point
$x\, \in\, S$ is the one in (\ref{eq4-1}), and the parabolic
weight of the subspace $W^j_{x}\,\subset\, W_{x}$ 
in (\ref{eq4-1}) is $(N-j)/N$.

Let $V_*$ be a parabolic vector bundle over $X$ satisfying
condition (1) and condition (2) in Theorem \ref{THE-1}. The
vector bundle underlying $V_*$ will be denoted by $V_0$.
Let $\widehat{V}$ be the $\Gamma$--linearized vector 
bundle over $Y$ corresponding to $V_*$. Then
\begin{equation}\label{e0-1}
H^0(Y,\, \widehat{V})^\Gamma\, =\, H^i(X,\, V_0)
\end{equation}
\cite[p. 310, (2.9)]{Bis1}.
Given any finite dimensional complex left
$\Gamma$--module $M$, there is a canonical
$\mathbb C$--linear isomorphism
\[
M\, \longrightarrow\,
\text{Hom}_{{\mathbb C}(\Gamma)}({\mathbb C}(\Gamma)\, ,M)
\, =\, \text{Hom}_{\mathbb C}({\mathbb C}(\Gamma)\, ,M)^\Gamma
\]
that sends any $v\, \in\, M$ to the homomorphism of
$\Gamma$--modules
$\rho_v\, :\, {\mathbb C}(\Gamma) \, \longrightarrow\, M$
uniquely determined by the condition that $\rho_v(\gamma)\, =\,
\gamma\cdot v$ for all $\gamma\, \in\, \Gamma$.
Using this canonical isomorphism, from (\ref{e0-1})
it follows that
\begin{equation}\label{e0-2}
H^0(Y,\, \widehat{V})\, =\,\text{Hom}_{\rm par}(W_*, V_*)\, .
\end{equation}

Since the bijective correspondence between the parabolic
vector bundles and the $\Gamma$--linearized vector bundles
is compatible with tensor product, dual and homomorphism,
the parabolic vector bundle corresponding to the
$\Gamma$--linearized vector bundle
${\widetilde R}^*_{r,d| \Gamma|}\bigotimes \wh F$ is
the parabolic tensor product
$(R'_{r,d,*})^*\bigotimes F_*$, where $(R'_{r,d,*})^*$
is the parabolic dual of $R'_{r,d,*}$.

Now the lemma follows by 
substituting the $\Gamma$--linearized vector bundle
${\widetilde R}^*_{r,d| \Gamma|}\bigotimes \wh F$ in place
of $\widehat{V}$ in (\ref{e0-2}).\hfill $\Box$
\medskip

Let 
\begin{equation}\label{R}
R_{r,d,*}\, :=\, W_*\otimes R'_{r,d,*}
\end{equation}
be the parabolic tensor product. In view of (\ref{eq4}) and
Lemma \ref{lem1}, we conclude that the parabolic vector bundle
$R_{r,d,*}$ constructed in (\ref{R})
satisfies the condition in Theorem \ref{THE-1}.\\

{\bf Step 3:} To complete the proof of Theorem \ref{THE-1}
we now consider the remaining case where
$(\text{genus}(X)\, , |S|)\, =\, (0\, ,1)$.
So $X\, =\, {\mathbb C}{\mathbb P}^1$ and $S \,=\,\{x\}$
is a singleton set.

Any vector bundle over ${\mathbb C}{\mathbb P}^1$ decomposes into
a direct sum of line bundles \cite[p. 126, Th\'eor\`eme 2.1]{Gr}.
Using this it follows immediately that a parabolic
vector bundle $E_*$ of rank $r$
on ${\mathbb C}{\mathbb P}^1$ with parabolic
structure over $x$ is parabolic semistable if and only if
the vector bundle $E_0$ underlying $E_*$ is a direct sum
$$
E_0\, =\, L^{\oplus r}\, ,
$$
and the quasiparabolic filtration is $0\, \subset\, (E_0)_x$.

Given $r$ and $d$ as in Theorem \ref{THE-1}, let
$$
d_0 \,=\, \lfloor d/r \rfloor
$$
be the integral part of $d/r$. Set
\begin{equation}\label{al.}
\alpha\, :=\, d/r-d_0\, \in\, [0\, ,1)\, .
\end{equation}
If $\alpha\, <\, (N-1)/N$, then set
$R_{r,d,*}$ to be the line bundle
${\mathcal O}_{{\mathbb C}{\mathbb P}^1}(d_0)$
equipped with parabolic weight $\alpha+1/N$ at $x$.

If $\alpha\, =\, (N-1)/N$, then set $R_{r,d,*}$ to be the
line bundle ${\mathcal O}_{{\mathbb C}{\mathbb P}^1}(d_0+1)$
with the trivial parabolic structure.

Take any parabolic vector bundle $E_*$ of rank $r$ and
parabolic degree $d$ on ${\mathbb C}{\mathbb P}^1$ with
parabolic structure on $x$. If $E_*$ is parabolic semistable,
then from the above observation that the underlying vector
bundle $E_0$ is a direct sum of $r$ copies of a line bundle
and the quasiparabolic flag of $E_*$ is trivial it follows
immediately that ${\rm Hom}_{\rm par}(R_{r,d,*}\, , E_*)
\, =\, 0$.

Now assume that
\begin{equation}\label{nz}
{\rm Hom}_{\rm par}(R_{r,d,*}\, , E_*)\,
\not=\, 0\, .
\end{equation}
Let
$$
E_0\, =\, \bigoplus_{i=1}^r L_i
$$
be a decomposition into a direct sum of line bundles
of the underlying vector bundle $E_0$. From (\ref{nz})
it follows immediately that $\text{degree}(L_i)$ is independent
of $i$. Therefore, $E_0$ is a direct sum of $r$ copies of a
line bundle. Again from (\ref{nz}) it follows that all
the parabolic weights are at most $\alpha$ defined
in (\ref{al.}). Therefore, the parabolic weight must be
$\alpha$ with multiplicity $r$. Consequently, $E_*$ is
parabolic semistable. This completes the
proof of Theorem \ref{THE-1}.

\end{document}